# A Statistical Measure of a Population's Propensity to Engage in Post-Purchase Online Word-of-Mouth

**Chrysanthos Dellarocas and Ritu Narayan**

*Abstract.* The emergence of online communities has enabled firms to monitor consumer-generated online word-of-mouth (WOM) in real-time by mining publicly available information from the Internet. A prerequisite for harnessing this new ability is the development of appropriate WOM metrics and the identification of relationships between such metrics and consumer behavior. Along these lines this paper introduces a metric of a purchasing population's propensity to rate a product online. Using data from a popular movie website we find that our metric exhibits several relationships that have been previously found to exist between aspects of a product and consumers' propensity to engage in offline WOM about it. Our study, thus, provides positive evidence for the validity of our metric as a proxy of a population's propensity to engage in post-purchase online WOM. Our results also suggest that the antecedents of offline and online WOM exhibit important similarities.

*Key words and phrases:* Word-of-mouth metrics, online communities, viral marketing, motion picture reviews.

## 1. INTRODUCTION

Among the many and varied channels through which a person may receive information, it is hard to imagine any that carry the credibility and importance of interpersonal communication, or word-of-mouth (WOM). Managers have always recognized WOM as an important driver of consumer behavior and have, thus, been interested in properly measuring and managing it. Nevertheless, the "perishable" nature of this important social phenomenon has, so far, limited the reliability of such efforts and resisted the development of a systematic discipline of WOM measurement.

Surveys have traditionally been the most popular method to measure WOM, largely because individuals can be asked directly about their communication habits; possible bias in measurement can then arise from the self-reporting of behavior. Several well-known studies, such as Bowman and Narayandas (2001), Brown and Reingen (1987), Reingen and Kernan (1986) and Richins (1983), base their analyses on proprietary surveys designed to test a specific hypothesis related to WOM.

The Internet has had a profound impact on both the reach as well as the measurability of WOM. The emergence of a diverse mosaic of Internet-mediated communities (product review sites, discussion groups, instant messaging chat rooms, web logs, etc.) has allowed individuals all over the world to easily share opinions on a variety of topics, including products, services, and even world events. Thanks to these systems, opinions of a single individual can instantly reach thousands, or even millions, of other consumers. This escalation in audience is altering the dynamics of many industries where WOM has traditionally played an important role. For example, the entertainment industry has found that the rapid spread of online WOM is shrinking the life cycles of its products (movies) and causing it to rethink its marketing strategies (Muñoz, 2003).

*Chrysanthos Dellarocas is Assistant Professor, Robert H. Smith School of Business, University of Maryland, College Park, Maryland 20742, USA e-mail: cdell@rhsmith.umd.edu. Ritu Narayan is a Doctoral Candidate, Robert H. Smith School of Business, University of Maryland, College Park, Maryland 20742, USA e-mail: rnarayan@rhsmith.umd.edu.*

At the same time, the Internet has made WOM instantly measurable: persistent traces of WOM (online product reviews, discussion group postings, web log entries, etc.) can be found in many publicly available Internet forums. This public data provides organizations with the ability to measure WOM *as it happens* by monitoring information available on the Internet.

A prerequisite for harnessing this newfound ability is the development of appropriate WOM metrics and the identification of relationships between such metrics and consumer behavior. In recent years a growing number of researchers and practitioners have been focusing on this important set of questions.

The majority of past research on this topic has focused on the use of online WOM as a revenue forecasting tool. Three metrics of online WOM have received particular attention in this context: volume, valence and dispersion. The theory behind measuring *volume*, or the number of online messages posted on a topic, is that the more consumers discuss a product, the higher the chance that other consumers will become aware of it. Liu (2004) found that the volume of messages posted on Internet message boards about upcoming and newly released movies was a good predictor of their box office success. The theory behind *valence*, or the fraction of positive and negative opinions in the mix of messages, is that, in addition to building awareness, WOM carries important information about a product's quality. Dellarocas, Awad and Zhang (2005) found that the valence of online ratings posted during a movie's opening weekend was the most important predictor of that movie's revenue trajectory in subsequent weeks. The theory behind measuring *dispersion*, or the spread of communication across communities, is that WOM spreads quickly *within* communities, but slowly *across* them (Granovetter, 1973). Godes and Mayzlin (2004) found that the dispersion of conversations about weekly TV shows across Internet communities has positive correlation with the evolution of viewership of these shows.

As firms become more sophisticated in managing WOM, they are increasingly engaging in proactive WOM-generation marketing campaigns (also known as *viral* marketing campaigns) whose objective is to induce adopters to "spread the word" about a product and, thus, to multiply product awareness and interest above and beyond what can be achieved through traditional marketing. An important prerequisite for the success of such campaigns is a better understanding of the factors that affect an individual's propensity to engage in WOM in relation to a stimulus she has received. Accordingly, the evaluation of such campaigns requires the development of metrics of a population's actual propensity to engage in product-related WOM.

In this paper we introduce one such metric. We call our new metric the *density* of online ratings; we define it as the ratio of the total number of people who posted online ratings for a product during a given time period over the number of people who bought that product during the same period. The density of online ratings, thus, represents a population-level estimate of the conditional probability that a person who has purchased a product will rate it online.

It is useful to think of density as somewhat analogous to a number of other metrics that are used to evaluate the effectiveness of marketing campaigns. For example, *click-through rates* are used to evaluate the effectiveness of online ads, while *Nielsen ratings* (viewership rates) are used to evaluate the effectiveness of television commercials. Density, similarly, offers marketing managers a tool to evaluate their WOM campaigns. Managers can evaluate if a WOM campaign has been successful by evaluating the propensity of people to talk about the product after the campaign.

In the rest of the paper we describe a study that provides initial evidence for the validity of our new metric. Using data collected from an online movie ratings community, we show that the density of online ratings exhibits a number of empirical relationships that have been previously shown to affect an individual's propensity to engage in offline WOM. Our results show that the density of online ratings can serve as a useful proxy of a purchasing population's propensity to engage in post-purchase WOM. They also suggest that the antecedents of offline and online WOM are similar.

Our choice of movies as the product category for observing online WOM offers some unique advantages. Unlike other product categories (such as books or music) for which online user ratings are easily accessible, but sales figures are not, production and weekly box office data are publicly available for most movies. Coupled with data on the average price of a movie ticket, this allows us to estimate the number of people who watched a particular movie during a specific week and relate this number to the corresponding volume of ratings posted online. [The motion picture industry is the focus of a rich academic literature (see, e.g., Neelamegham and Chintagunta, 1999; Simonoff and Sparrow, 2000; Eliashberg and Shugan, 1997; and references contained in Dellarocas, Awad and Zhang, 2005), primarily focused on understanding the factors that affect a movie's box office success. Our work relates to that literature but addresses a question (what factors make moviegoers more or less eager to post an online movie review?) that has not been studied by any previous work.]

## 2. DATA

Data for this study were collected from Yahoo! Movies (movies.yahoo.com) and BoxOfficeMojo (www.boxofficemojo.com). From Yahoo! Movies, we collected the names of all movies released during 2002. For the purpose of our analysis, we excluded titles that were (a) not released in the United States, (b) not a theatrical release (e.g., DVD releases), or (c) not released nationwide. For each of the remaining titles we collected detailed ratings information, including all professional critic reviews (text and letter ratings, which we converted to a number between 1 and 5) and all user reviews (date and time of review, user id, review text, integer ratings between 1 and 5).

We used BoxOfficeMojo to obtain weekly box office, budget and marketing expenses data. This information was missing for several movies. Our final dataset consists of 104 movies, 1392 critic reviews (an average of 13 reviews per movie) and 63,889 user reviews from 46,294 individual users (an average of 614 reviews per movie and 1.4 reviews per user). Table 1 provides some key summary statistics.

Table 1 shows that the mean and standard deviation of ratings posted by users and professional critics across all movies of our dataset are very similar. This statistical fact hides an important difference in the rating behavior of users and critics. Recall that movie scores are integers between 1 and 5. Figure 1 plots the relative use of these scores by users and critics. Interestingly, whereas critics seem to be rating movies on a (slightly upwardly biased) curve, the majority of user ratings lie at the two extremes of the ratings scale, with a strong emphasis on the positive end: almost half of all posted ratings are equal to the highest possible rating, 18% of ratings are equal to the lowest possible rating, and only about 30% are intermediate values. The preponderance of extreme reviews is consistent with past research on word-of-mouth that finds that people are more likely to engage in interpersonal communication when they have very positive and very negative experiences (Anderson, 1998).

We were able to collect partial rater demographic data by mining the user profiles that are associated with the raters' Yahoo! ids. Eighty-five percent of raters in our sample listed their gender and 34% their age. From that information, we constructed a noisy estimate of the demographic profile of the Yahoo! Movies rater population. We found the demo-graphic breakdown of online raters to be substantially skewed relative to the profile of US moviegoers (as reported by the Newspaper Association of America). Most notably, a disproportionately high percentage of online ratings were provided by males (an estimated 74% of raters, compared with only 49% of US moviegoers) and by people between 18 and 29 years of age (an estimated 58% of raters, compared with only 35% of US moviegoers).

TABLE 1
*Key summary statistics of our dataset*

| Variable | Mean | Std. dev. | Min | Max |
|---|---|---|---|---|
| Marketing budget (in millions) | 24.09 | 11.67 | 1 | 50 |
| Screens in opening week | 1964.69 | 1101.88 | 4 | 3876 |
| Weekly box office revenue (in millions) | 11.80 | 17.43 | 0.014 | 151.62 |
| Weekly volume of user ratings | 124.97 | 282.22 | 1 | 3802 |
| Weekly ratings density (ratings per million viewers) | 79 | 118.80 | 4 | 1738 |
| Weekly average of user ratings (range 1–5) | 3.44 | 0.69 | 1 | 5 |
| Average of critics' ratings (range 1–5) | 3.35 | 0.68 | 1.29 | 4.28 |
| Total number of movies | 104 | | | |
| Total number of user ratings | 63,889 | | | |
| Total number of critic ratings | 1,392 | | | |
| Total number of unique users | 46,294 | | | |

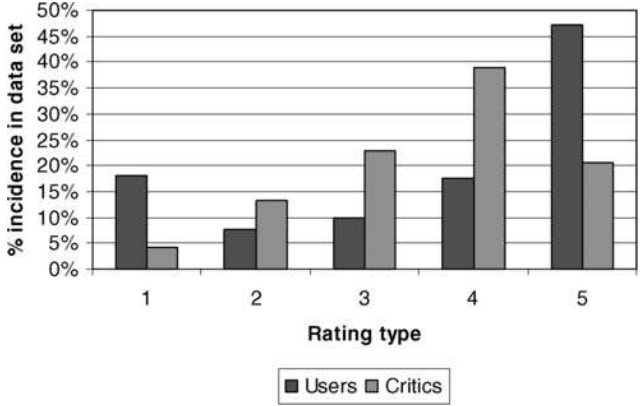

Fig. 1. *Relative use of different scores by users and critics.*

## 3. THE WOM DENSITY METRIC

Our interest in this study focuses on measuring a moviegoer's propensity to post an online rating for a movie she has recently watched. Mathematically, consumer $i$'s propensity to rate movie $j$ corresponds to the conditional probability

$$Pr[i \text{ rates } j | i \text{ watched } j]$$
$$= \frac{Pr[i \text{ watched } j | i \text{ rates } j] Pr[i \text{ rates } j]}{Pr[i \text{ watched } j]}$$
$$= \frac{Pr[i \text{ rates } j]}{Pr[i \text{ watched } j]}$$

where we assume that $Pr[i \text{ watched } j | i \text{ rates } j] = 1$, that is, that people only post ratings for movies they have already watched. Even though we do not have data on individual raters' movie attendance, a population-level estimate of the above conditional probability is simply the *density* of online ratings, defined as

$D_{jt} =$ number of people who posted ratings for movie $j$ during period $t$

$\cdot$ (number of people who watched movie $j$ during period $t)^{-1}$

$=$ number of ratings posted for movie $j$ by unique users during period $t$

$\cdot$ (average ticket price)

$\cdot$ (box office revenues of movie $j$ during period $t)^{-1}$.

The above definition of density assumes that:

1. People post at most one rating for a given movie.
2. People only watch a movie once during the period of observation.
3. People who watch a movie during period $t$ post online ratings within the same period (or not at all).

The first assumption can be easily satisfied if we only "count" one rating per movie and user. Under the assumption that most users post all their ratings under the same online pseudonym, this is easy to do since Yahoo! movie ratings have an associated Yahoo! id. The second assumption restricts our density metric to durable goods, that is, goods that are purchased only once, or very infrequently. Movies fall into this category since it seems reasonable to assume that most people watch a movie only once within our window of observation (five weeks). Finally, our data provides strong support for the third assumption when $t$ is considered on a weekly level. Figure 2 plots the daily box office revenues and corresponding daily volume of ratings for "Spider-Man." Observe that the volume of ratings closely follows the box office peaks and valleys, suggesting that most moviegoers posted ratings within a week of watching the movie. Most movies in our data set exhibit similar patterns. To establish this formally, we calculated the correlation between the weekly volume of ratings and the lagged weekly box office revenues. Correlation was highest between volume and box office revenues of the same week and monotonically declined as the lag increased.

Despite the high correlation between box office revenues and ratings volume, the density of individual movies exhibits substantial variance. Figure 3 plots the empirical cumulative probability distribution of first-week's ratings density for movies in our dataset, whereas Table 2 lists the movies with the highest and lowest densities. First-week density ranges between 14 and 701 ratings per million viewers; the highest density is, thus, about 50 times higher than the lowest density. The highest first-week density in our sample corresponds to "Swept

Away," a drama/romance movie starring Madonna that received almost universally bad reviews and flopped in the box office. It was followed by "Solaris," a highly praised US remake of a cult Russian science fiction movie, and "Frida," a well-received movie about the life of Mexican artist Frida Kahlo produced by and starring Salma Hayek, that was released as a "sleeper." (Most movies are distributed using one of two strategies. Wide-release or "blockbuster" movies are released on a large number of screens and are accompanied by extensive marketing campaigns. Box office revenues for such movies are typically highest during the first week and steadily decline afterward. Narrow-release or "sleeper" movies, on the other hand, are initially released on a small number of screens, have more modest marketing budgets and rely on word-of-mouth for revenue growth.) Most movies at the low end of the density scale were children's movies.

Managers planning or evaluating viral marketing campaigns are interested in better understanding the factors that explain such large variance in consumers' propensity to discuss movies online. The remainder of the paper provides some initial answers to this question.

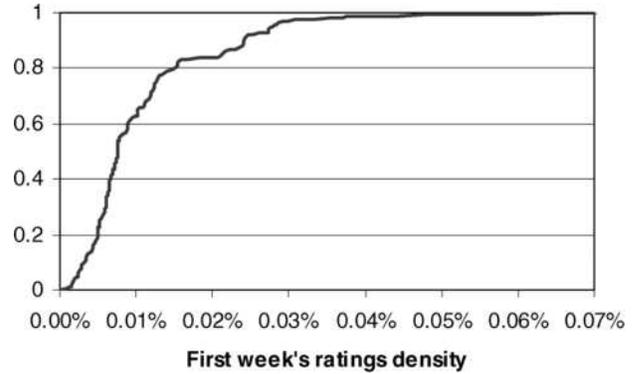

FIG. 3. *Empirical cumulative distribution of first-week's ratings density in our dataset.*

## 4. PROPERTIES OF THE DENSITY METRIC

### 4.1 Hypotheses

Even though the systematic measurement of WOM using Internet data is a recent development, the study of the antecedents of (offline) WOM using traditional methodologies is much older. A long stream of studies has established a number of factors that affect an individual's propensity to engage in post-purchase WOM. We provide evidence for the validity of the density metric as a proxy of a purchasing population's propensity to engage in online WOM, by studying to what extent it exhibits relationships consistent with those previously established in the study of offline WOM.

A number of studies have shown that consumers who experience extreme satisfaction or dissatisfaction with a product exhibit higher propensity to en-gage in WOM (see, e.g., Anderson, 1998). We expect to find a similar relationship online.

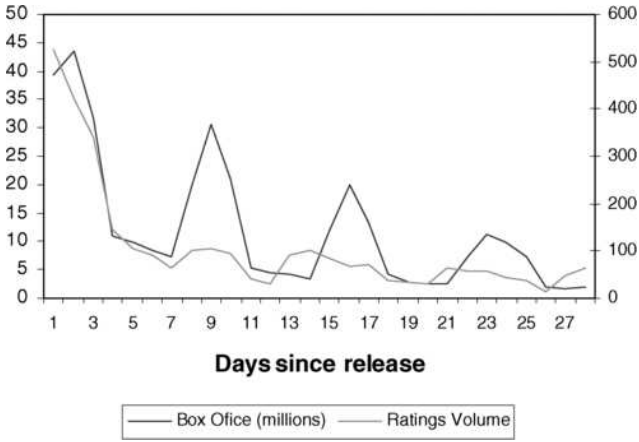

FIG. 2. *Daily box office revenues (left axis) and corresponding daily volume of ratings (right axis) for "Spider-Man."*

TABLE 2
*Movies with highest and lowest ratings density*

| Top 5 | | | Bottom 5 | | |
|---|---|---|---|---|---|
| Title | Genre | Density[a] | Title | Genre | Density[a] |
| Swept Away | Romance | 701 | Spirit: Stallion of the Cimarron | Children | 23 |
| Solaris | Sci-Fi | 380 | Fear Dot Com | Thriller | 19 |
| Frida | Drama | 371 | Spy Kids 2 | Children | 18 |
| Confessions of a Dangerous Mind | Comedy Drama | 292 | Stuart Little 2 | Children | 16 |
| About Schmidt | Drama | 284 | Return to Never Land | Children | 14 |

[a]Yahoo! ratings per million viewers.

H1: *The density of online ratings will be higher for movies that are perceived by consumers to be exceptionally good or exceptionally bad.*

Marketing campaigns have been found to stimulate word-of-mouth. A number of past studies have provided evidence for this theory, finding that repetitive advertising (Bayus, 1985) and unusual advertising (King and Tinkham, 1990) increase the propensity of consumers to engage in WOM. If one were to expect a similar pattern in online ratings, the following hypothesis should hold:

H2: *A movie's density of online ratings is positively related to that movie's marketing effort.*

Several authors have argued that consumers have a higher propensity to share their experiences with others when these make them look intelligent and savvy (Dichter, 1966; Sundaram, Mitra and Webster, 1998; Hennig-Thurau et al., 2004). This would predict a higher tendency to post online reviews about more eclectic, less widely released movies.

H3: *A movie's density of online ratings is negatively related to the availability of that movie.*

Finally, past research has argued that people often engage in WOM to reduce doubts about the purchase of controversial products (Engel, Blackwell and Miniard, 1995). According to this *dissonance reduction* theory one expects that the propensity to engage in WOM will be higher when there is more public disagreement surrounding a product. This leads to the following hypothesis:

H4: *A movie's density of online ratings is positively related to the amount of public disagreement about that movie's quality.*

### 4.2 Independent variables

We use a movie's cumulative marketing budget ($MKT$) as a proxy of the marketing effort related to the movie. The number of screens on which the movie is exhibited in week $t$ ($SCR_t$) serves as a proxy of a movie's availability. The arithmetic mean of user ratings posted for that movie during week $t$ ($AVG_t$) serves as a proxy of perceived quality. The amount of public controversy that surrounds the perceived quality of a movie is proxied by the standard deviation of professional critic ratings ($CRSTD$). [In theory, the standard deviation of user ratings would have been an even better measure of disagreement. However, the distribution of user ratings for most movies is severely bimodal (Figure 1). It is well known that the standard deviation of such distributions is not a good measure of dispersion and exhibits high correlation to the mean. On the other hand, the distribution of critic reviews for most movies is unimodal. Its standard deviation is, thus, a better measure of disagreement.]

Dummy variables are used to control for movie genre ($ROMANCE$, $SCIFI$, $DRAMA$, $THRILLER$, $COMEDY$, $ACTION$, $KIDS$). Since a movie can belong to multiple genres, in our coding scheme for the genre we assign a score of $1/k$ for each genre for a movie that falls in $k$ genres. Finally, we use the log transformation of the number of weeks since release ($LWEEK$) to control for the decay in WOM about a movie in later weeks. The use of the log transformation was empirically motivated by observing that, whereas the density of most movies in our dataset drops substantially between the first and second week of release, the rate of decay slows down after the second week.

### 4.3 Model specification and results

Since $D_{jt}$ represents a proportion, we apply a logit transformation to linearize its relationship with independent variables and improve model fit (Greene, 2003). In the rest of the study we will use $LD_{jt} = \log(D_{jt}/(1-D_{jt}))$ as our dependent variable, with $t = 1, 2, \ldots$ referring to weeks elapsed since the movie's national release. We restrict our attention to weeks 1–5 since WOM effects diminish after this time period.

Our specification is a combined model of (the logit transformation of) weekly density:

$$\begin{aligned}
LD_{jt} = {} & \beta_0 + \beta_1 MKT_j + \beta_2 SCR_{jt} + \beta_3 AVG_{jt} \\
& + \beta_4 (AVG_{jt})^2 + \beta_5 CRSTD_j \\
& + \beta_6 ROMANCE_j + \beta_7 THRILLER_j \\
& + \beta_8 DRAMA_j + \beta_9 COMEDY_j \\
& + \beta_{10} SCIFI_j + \beta_{11} ACTION_j \\
& + \beta_{12} KIDS_j + \beta_{13} LWEEK_{jt} + \varepsilon_{jt}.
\end{aligned} \quad (1)$$

To test for the presence of a U-shaped relationship between perceived quality and ratings density (Hypothesis 1) we included both first- and second-degree terms of average user ratings ($AVG_{jt}$). To alleviate potential multicollinearity problems between the first- and second-degree terms, we used mean-centered values for the user rating variables by subtracting their mean from each of the observations.

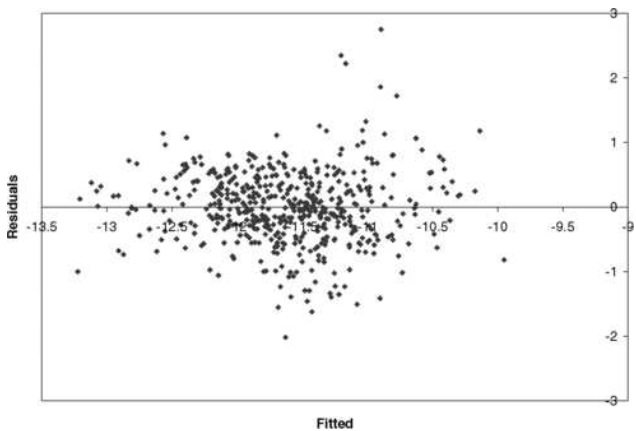

FIG. 4. *Residuals vs. fitted values after OLS regression with logit transformation.*

Heteroskedasticity is, almost by definition, present in our model. Let $N_{jt}$ denote the number of people who watched movie $j$ during week $t$. ($N_{jt}$ is equal to the corresponding box office revenues divided by the average ticket price.) Our dependent variable, $D_{jt}$, which is meant to be an estimate of the conditional probability $p_{jt} = Pr[$representative user rates movie $j$ during week $t|$user watched $j]$, can equivalently be thought of as the sample mean of $N_{jt}$ independent binary random variables (representing each individual moviegoer's decision to rate or not rate a movie) drawn from Bernoulli distributions with probabilities of success $p_{jt}$. According to this perspective, the variance of $D_{jt}$ is proportional to $1/N_{jt}$, that is, inversely proportional to movie $j$'s box office revenues during week $t$. Since weekly box office revenues vary substantially among movies in our sample, heteroskedasticity is likely to be significant.

Our data confirms this hypothesis. First, we fitted our model using ordinary least squares regression. Inspection of the plot of residuals versus fitted values (Figure 4) indicates presence of heteroskedastic errors. We subsequently used the Breusch–Pagan test (Breusch and Pagan, 1979) to formally test for the presence of heteroskedasticity and found positive evidence. We address this issue by performing weighted least squares regression (Greene, 2003); consistent with the above discussion the weights we assign to each observation ($D_{jt}$) are proportional to the corresponding box office revenues of movie $j$ during week $t$.

Table 3 summarizes the results of our regression analysis. Model fit is reasonably good ($Adj. R^2 = 0.55$) and all our critical variables are statistically significant. Next, we describe our key findings:

1. The coefficient of $MKT_j$ is positive. This indicates that marketing expenditures have a positive relationship with ratings density.
2. The coefficient of $SCR_{jt}$ is negative. This finding indicates a negative relationship between a movie's availability and the density of its online ratings.
3. The coefficient of the second-order term $(AVG_{jt})^2$ is positive, whereas the coefficient of the corresponding first-order term is not significant. We, thus, find support for the presence of a U-shaped relationship between perceived movie quality, as expressed by the average user ratings, and observed ratings density (Figure 5).
4. The coefficient of $CRSTD_j$ is positive, indicating that higher disagreement among critic reviews relates to higher ratings density for that movie.

Regarding the fixed effects of movie genre, after all above factors have been taken into account, Table 3 shows that science-fiction movies have a significant positive effect in moviegoers' propensity to post online ratings, whereas children's and comedy movies have a significant negative effect. Drama and action movies also have a negative effect. The coefficients for movies in thriller and romance genre are insignificant. Given the demographics of our rater population (disproportionate fraction of young male adults), these results are in line with common sense and also help explain some of the findings of Table 2. Finally, as expected, the coefficient of $LWEEK$ is negative and significant since the propensity of

TABLE 3
*WLS regression results*

| Variable | Coeff. | Std. err. | $t$-value | $P > |t|$ | |
|---|---|---|---|---|---|
| MKT | 0.02 | 0.003 | 7.30 | <0.001 | c |
| SCR | −0.0005 | 0.00003 | −12.43 | <0.001 | c |
| AVG | 0.07 | 0.04 | 1.78 | 0.076 | |
| AVG$^2$ | 0.11 | 0.03 | 3.05 | 0.002 | b |
| CRSTD | 0.33 | 0.07 | 5.71 | <0.001 | c |
| SCIFI | 0.54 | 0.06 | 8.76 | <0.001 | c |
| KIDS | −0.82 | 0.07 | −11.99 | <0.001 | c |
| DRAMA | −0.15 | 0.06 | −2.58 | 0.010 | a |
| COMEDY | −0.21 | 0.06 | −3.58 | <0.001 | c |
| ROMANCE | −0.11 | 0.08 | −1.38 | 0.167 | |
| ACTION | −0.15 | 0.06 | −2.68 | 0.008 | b |
| THRILLER | 0.09 | 0.06 | 1.62 | 0.107 | |
| LWEEK | −0.55 | 0.04 | −13.71 | <0.001 | c |
| _CONS | −11.14 | 0.16 | −70.33 | <0.001 | c |

*Significance codes: 0[c], 0.001[b], 0.01[a].

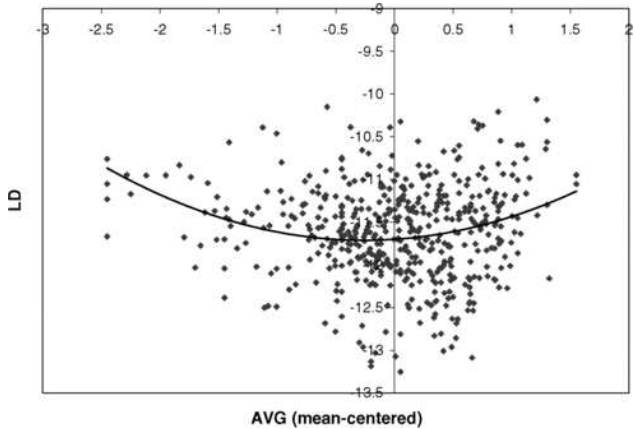

FIG. 5. *Plot of U-shaped relationship between fitted values of logit transformation of density and average weekly user rating.*

moviegoers to talk about a movie tends to decay as the movie gets "old."

## 5. CONCLUDING REMARKS

Our findings show that the density of online ratings exhibits several relationships that have been previously found to exist between aspects of a product and the consumers' propensity to engage in offline WOM related to that product. Our study, thus, provides positive evidence for the usefulness of the density metric as a proxy of a purchasing population's propensity to engage in post-purchase online WOM. Our results also suggest that the antecedents of offline and online WOM are similar. If true, such a result is important because it means that several insights obtained through decades of research on offline WOM can apply to the online domain. These insights, together with the new ability to measure aspects of WOM in real-time by mining publicly available data from Internet communities, can lead to substantial advances in the ability of organizations to manage WOM. The science of statistics has an important role to play in this nascent endeavor in terms of discovering and evaluating new metrics and proposing appropriate techniques for analyzing the huge volume of online ratings data.

Our results have interesting implications for viral marketing campaigns, that is, campaigns that aim to generate online WOM. Extreme (dis)satisfaction, controversy, advertising and product exclusivity all seem to correlate with higher propensity to discuss a product online. The impact of exclusivity in our dataset is important because, in most cases, higher marketing spending is positively correlated with broader product availability. Our initial findings, thus, suggest that product marketing campaigns that aim to maximize word-of-mouth might want to strike a balance between advertising spending and maintaining an image of exclusivity. Since our study does not establish causality, further research is needed to prove or disprove this hypothesis.

We conclude by mentioning some of our study's limitations and associated directions for further research. First, sampling biases in our data could have stemmed from Yahoo! Movies not representing all websites where movie reviews are posted, plus the fact that the 104 movies in our dataset only represent a subset of movies released in 2002. Further, our study simply identifies associations between the density of online ratings and other attributes of a movie. Establishing causality would require a substantially more complex modeling approach. Another limitation of our dataset is that we do not have information related to the set of movies that individual raters watched (but never rated). This was not a handicap in estimating the population's average propensity to rate but would pose important challenges if we wanted to delve deeper into the population and study, say, the presence of segments whose rating behavior follows different patterns. Finally, our study only looks at one product category. Our objective is to apply the techniques explored in this paper to other industries, in order to better understand what makes people engage in online WOM.

Additional insights can be gained by modeling the dynamics of a population's propensity to rate movies online. Such analyses would benefit from the use of advanced statistical techniques, such as functional data analysis (see, e.g., Wu, 2005). A better understanding of the dynamics of density would allow us to distinguish between products that generate short-lived consumer WOM from those that induce a more sustained consumer response as well as to characterize the properties of online WOM in early versus later stages of a new product's release. Similar analyses would allow us to study aspects unique to online WOM, such as the impact of the presence of prior ratings on the arrival rate and valence of subsequent ratings.

## REFERENCES

ANDERSON, E. W. (1998). Customer satisfaction and word of mouth. *J. Service Research* **1** 5–17.

BAYUS, B. L. (1985). Word-of-mouth: The indirect effects of marketing efforts. *J. Advertising Research* **25** 31–39.

Bowman, D. and Narayandas, D. (2001). Managing customer-initiated contacts with manufacturers: The impact on share of category requirements and word-of-mouth behavior. *J. Marketing Research* **38** 281–297.

Breusch, T. and Pagan, A. (1979). A simple test for heteroscedasticity and random coefficient variation. *Econometrica* **47** 1287–1294. MR0545960

Brown, J. J. and Reingen, P. (1987). Social ties and word-of-mouth referral behavior. *J. Consumer Research* **14** 350–362.

Dellarocas, C., Awad, N. and Zhang, M. (2005). Using online ratings as a proxy of word-of-mouth in motion picture revenue forecasting. Working paper, Smith School of Business, Univ. Maryland.

Dichter, E. (1966). How word-of-mouth advertising works. *Harvard Business Review* **44**(6) 147–160.

Eliashberg, J. and Shugan, S. M. (1997). Film critics: Influencers or predictors? *J. Marketing* **61**(2) 68–78.

Engel, J. F., Blackwell, R. D. and Miniard, P. W. (1995). *Consumer Behavior*, 8th ed. Dryden, Fort Worth.

Godes, D. and Mayzlin D. (2004). Using online conversations to study word of mouth communication. *Marketing Science* **23** 545–560.

Granovetter, M. (1973). The strength of weak ties. *American J. Sociology* **78** 1360–1380.

Greene, W. (2003). *Econometric Analysis*, 5th ed. Prentice Hall, Upper Saddle River, NJ.

Hennig-Thurau, T., Gwinner, K. P., Walsh, G. and Gremler, D. D. (2004). Electronic word-of-mouth via consumer-opinion platforms: What motivates consumers to articulate themselves on the Internet? *J. Interactive Marketing* **18** 38–52.

King, K. W. and Tinkham, S. F. (1990). The learning and retention of outdoor advertising. *J. Advertising Research* **29** 47–51.

Liu, Y. (2004). Word-of-mouth for movies: Its dynamics and impact on box office receipts. Working paper.

Muñoz, L. (2003). High-tech word of mouth maims movies in a flash. *Los Angeles Times* August 17.

Neelamegham, R. and Chintagunta, P. (1999). A Bayesian model to forecast new product performance in domestic and international markets. *Marketing Science* **18** 115–136.

Reingen, P. and Kernan, J. (1986). Analysis of referral networks in marketing: Methods and illustration. *J. Marketing Research* **23** 370–378.

Richins, M. L. (1983). Negative word-of-mouth by dissatisfied consumers: A pilot study. *J. Marketing* **47**(1) 68–78.

Simonoff, J. S. and Sparrow, I. R. (2000). Predicting movie grosses: Winners and losers, blockbusters and sleepers. *Chance* **13**(3) 15–24. MR1854211

Sundaram, D. S., Mitra, K. and Webster, C. (1998). Word-of-mouth communications: A motivational analysis. *Adv. in Consumer Research* **25** 527–531.

Wu, O. (2005). Dynamics of online movie ratings. Term paper, Research Interactive Team, VIGRE program, Univ. Maryland.